\documentclass{amsart}
\usepackage{enumerate,amssymb,amsmath,graphics,epsfig}
\usepackage{color}
\usepackage{xspace}
\usepackage{mathrsfs}
\usepackage{arydshln}
\usepackage[labelfont=rm,font=small,labelformat=simple]{subcaption}

\usepackage{parskip}
\usepackage{pslatex}
\usepackage{amsfonts}
\usepackage{tikz}
\usepackage{lmodern}
\usepackage{graphicx}
\usepackage{bm}
\usepackage{ulem}
\usepackage{setspace}
\usepackage{hyperref}
\usepackage{caption}
\usepackage{subcaption}
\usepackage{marginnote}
\usepackage[foot]{amsaddr}
\newcommand\gA{\underline{\mathbb{A}}}
\newcommand\gB{\underline{\mathbb{B}}}
\newcommand\gC{\underline{\mathbb{C}}}
\newcommand\gD{\underline{\mathbb{D}}}
\newcommand\nab{\underline{\bm{\nabla}}}
\theoremstyle{plain}

\theoremstyle{remark}

\begin{document}
%%%%%%%%%%%%%%%%%%%%%%%%%%%%%%%%%%%%%%%%%%%%%%%%%%%%%%%%%%%%%%%%%%%%%%%%%%%%%%%%
\title[]
{On the spectrum of the finite element approximation of a three field
formulation for linear elasticity}
\author{Linda Alzaben~$^{\ast}$}
\author{Daniele Boffi~$^\dagger$}
\address[$\ast$]{King Abdullah University of Science and Technology (KAUST), Saudi Arabia}
\address[$\dagger$]{King Abdullah University of Science and Technology (KAUST), Saudi Arabia and Universit\`a degli Studi di Pavia, Italy}
\email[$\ast$]{linda.alzaben@kaust.edu.sa}
\email[$\dagger$]{daniele.boffi@kaust.edu.sa}
%%%%%%%%%%%%%%%%%%%%%%%%%%%%%%%%%%%%%%%%%%%%%%%%%%%%%%%%%%%%%%%%%%%%%%%%%%%%%%%%
\begin{abstract}
We continue the investigation on the spectrum of operators arising from the discretization of partial differential equations.
In this paper we consider a three field formulation recently introduced for the finite element
least-squares approximation of linear elasticity.
We discuss in particular the distribution of the discrete eigenvalues in the complex plane and how they approximate the positive real eigenvalues of the continuous problem.
The dependence of the spectrum on the Lam\'e parameters is considered as well and its behavior when approaching the incompressible limit.
\end{abstract}
%%%%%%%%%%%%%%%%%%%%%%%%%%%%%%%%%%%%%%%%%%%%%%%%%%%%%%%%%%%%%%%%%%%%%%%%%%%%%%%%
\maketitle
\section{Introduction}
%%%%%%%%%%%%%%%%%%%%%%%%%%%%%%%%%%%%%%%%%%%%%%%%%%%%%%%%%%%%%%%%%%%%%%%%%%%%%%%%
Following~\cite{LZ_DB} we study the behavior of the spectrum of operators arising from the discretization of partial differential operators.
We consider a three field formulation of linear elasticity introduced in~\cite{BC}, based
on a least squares approach, and its finite element approximation.

While the continuous problem is self-adjoint, and its eigenvalues are positive real numbers, the approximation is based on a formulation that is not symmetric.
The asymmetry comes, as for other least squares formulations~\cite{FB_eig,FB_LE,DPG}, from the right hand side of the (singular) generalized algebraic eigenvalue problem.
As a consequence, we might expect that the eigenvalues are spread in the entire complex plane.
On the other hand, the uniform convergence of the discrete solution operator towards the continuous one implies that the eigenvalues are well approximated and that there are no spurious modes.
The consequences of this implication are numerically investigated; it turns out that, as expected, given a positive number $R$, for $h$ small enough, the number of discrete eigenvalues with modulus smaller than $R$ coincides with the number of continuous eigenvalues, when counted with their multiplicities.
There may be several complex eigenvalues with modulus larger than $R$ and this is clearly observed in our numerical tests.

We test that the convergence rate of the eigenvalues corresponds to the theoretical one and we study the behavior of the spectrum when the properties of the material approach the incompressible limit.
It turns out that, similarly to what we have observed in~\cite{LZ_DB}, the discrete spectrum
becomes more chaotic as the material is less compressible.
This remark could be interpreted, for instance, as a warning when this scheme is used for the
approximation of transient problem, whose behavior is clearly related to the properties of the discrete spectrum of the corresponding stationary problem.

In Section~\ref{sc:Contiuouse_porblem}, we recall the strong form of the linear elasticity problem and the energy functional associated with the three-field formulation of~\cite{BC}.
In Section~\ref{sec:three_field_eigenvalue_formulation}, we describe the variational formulation associated with the three-field formulation and in Section~\ref{sec:FiniteElementApproximation} its finite element approximation.
Finally, our main results are reported in Section~\ref{se:num}, where several numerical tests are performed.
%
%%%%%%%%%%%%%%%%%%%%%%%%%%%%%%%%%%%%%%%%%%%%%%%%%%%%%%%%%%%%%%%%%%%%%%%%%%%%%%%%
\section{Elasticity equations}
\label{sc:Contiuouse_porblem}
%%%%%%%%%%%%%%%%%%%%%%%%%%%%%%%%%%%%%%%%%%%%%%%%%%%%%%%%%%%%%%%%%%%%%%%%%%%%%%%%
%
In this section, we outline the stress-displacement formulation of linear elasticity and give a brief description on the problem.
Consider a bounded polygonal domain $\Omega$ in $\mathbb{R}^d,(d = 2,3)$ with boundary $\Gamma=\partial \Omega$.
The boundary is partitioned into two distinct open subsets $\Gamma_D$ and $\Gamma_N$ with $\partial\Omega = \overline{\Gamma}_D \cup\overline{\Gamma}_N$.
The strong form of the boundary-value problem is written as a first order system: find a displacement vector-field ${\bf{u}}_{(d\times 1)}$ and a symmetric stress tensor $\underline{\bm{\sigma}}_{(d\times d)}$ such that
\begin{equation}\label{eq:elasticity_First_order_system}
\begin{cases}
\mathcal{A}\underline{\bm{\sigma}} - 
\underline{\bm{\epsilon}}(\mathbf{u}) 
=
\bf{0} &\quad \text{in} \quad \Omega \\
\bf{div}\underline{\bm{\sigma}} =
- \mathbf{f}
&\quad \text{in} \quad \Omega \\
\end{cases}
\end{equation}
with homogeneous boundary conditions
\begin{equation} 
\begin{cases}
\mathbf{u} = 
\bf{0} \quad &\quad \text{on} \quad \Gamma_{D} \\
\mathbf{n} \cdot \underline{\bm{\sigma}} = 
\bf{0} &\quad \text{on} \quad{\Gamma}_{N}
\end{cases}
\end{equation}
where $\bf{n}$ is the outward unit normal vector to the boundary, $\underline{\bm{\epsilon}}(\bf{v})$ is the strain tensor, $\bf{f}$ the external force, and $\mathcal{A}$ is the compliance tensor given by
\begin{equation}\label{eq:compliance_tensor}
\mathcal{A}\underline{\bm{\tau}} =
\frac{1}{2\mu}\bigg(\underline{\bm{\tau}} -
\frac{\lambda}{d\lambda + 
2\mu}\text{tr}(\underline{\bm\tau})\underline{\mathbf{I}}\bigg).
\end{equation}
We denote by $\text{tr}(\cdot)$ and $\underline{\mathbf{I}}$ the trace of a matrix and the $d\times d$ identity matrix, respectively.
For isotropic elastic material, the compliance tensor is a fourth order tensor, which depends on the Lam\'e constants $\mu$ and $\lambda$.
Those constants are such that $\mu \in [\mu_1, \mu_2] $ with $0<\mu_1<\mu_2<\infty$ and $0<\lambda\leq \infty$.
When the value of $\lambda$ is large or infinity, the material is then said to be nearly incompressible or incompressible, respectively.
In this case, the compliance tensor is not invertible. 
The strain tensor $\underline{\bm{\epsilon}}({\mathbf{v}})$, is the symmetric
gradient of $\mathbf{v}$, that is, the symmetric part of $\nab \bf{v}$ given by
\begin{equation}
\underline{\bm{\epsilon}}({\mathbf{v}}) = 
\frac{1}{2}(\nab\mathbf{v} + 
(\nab\mathbf{v})^t).
\end{equation}
An equivalent representation of the symmetric gradient can be obtained by taking the difference between $\nab \bf{v}$ and its skew-symmetric part as follows:
\begin{equation}
\underline{\bm{\epsilon}}(\textbf{v}) = 
\nab\textbf{v}-
\underbrace{\frac{1}{2}(\nab\textbf{v}-(\nab\textbf{v})^t}_{\underline{\bm{\chi}}\bm{\omega} })= 
\nab\textbf{v} -(-1)^d \underline{\bm{\chi}}\bm{\omega},
\end{equation}
where $\bm{\omega}=\frac{1}{2}\nab \times \bf{v}$ denotes the vorticity, which is a scalar for $d=2$ and a vector for $d=3$ (for definition see~\cite{BC} Sec. 2).
Moreover, $\underline{\bm{\chi}}$ is the matrix given by
\begin{equation}
\underline{\bm{\chi}}=
\begin{cases}
\begin{pmatrix} 
0 & -1 \\ 
1 &  0
\end{pmatrix} & \text{if } d=2\\
(\underline{\bm{\chi}}_1,\underline{\bm{\chi}}_2,\underline{\bm{\chi}}_3) & \text{if } d=3
\end{cases}
\end{equation}
with
\begin{equation}
\underline{\bm{\chi}}_1 = 
\begin{pmatrix} 
0&0&0 \\ 
0&0&-1\\
0&1&0
\end{pmatrix},~
\underline{\bm{\chi}}_2 = 
\begin{pmatrix} 
0&0&1 \\ 
0&0&0\\
-1&0&0
\end{pmatrix}, \text{ and }
\underline{\bm{\chi}}_3 = 
\begin{pmatrix} 
0&-1&0\\ 
1&0&0\\
0&0&0
\end{pmatrix}.
\end{equation}
Considering the vorticity is an independent variable, the formulation in 
\eqref{eq:elasticity_First_order_system}, with the same boundary conditions, can be rewritten as follows:
\begin{equation}\label{eq:elasticity_with_vorticity}
\begin{cases}
\mathcal{A}\underline{\bm{\sigma}} - 
\nab\textbf{u}+(-1)^d\underline{\bm{\chi}}\bm{\omega}
=
\bf{0} &\quad \text{in} \quad \Omega \\
\bf{div}\underline{\bm{\sigma}} =
- \mathbf{f}
&\quad \text{in} \quad \Omega \\
\text{as}\underline{\bm{\sigma}} =
\mathbf{0}
&\quad \text{in} \quad \Omega 
\end{cases}
\end{equation}
The last equation in System \eqref{eq:elasticity_with_vorticity} ensures the symmetry of the stress tensor, since adding the vorticity in the first equation no longer implies its symmetry.
The value $\text{as} \underline{\bm{\sigma}}$ denotes the skew-symmetric part of $\underline{\bm{\sigma}}$.
Formulation \ref{eq:elasticity_with_vorticity} has the important characteristic of being robust in the incompressible limit as proven in~\cite{Cai}.

We define the solution spaces as follows: let
\begin{equation}
\underline{\textbf{X}} =
\begin{cases}
\textbf{H}(\textbf{div};\Omega)^d &\quad \text{if} \quad \Gamma_N\neq \emptyset\\
\{\underline{\bm{\tau}}\in\textbf{H}(\textbf{div};\Omega)^d:
\int_{\Omega}\text{tr}(\underline{\bm{\tau}})d\textbf{x} = 
0\} &\quad \text{if} \quad\Gamma_N= \emptyset\\
\end{cases}
\end{equation}
with the subspace
\begin{equation}
\underline{\textbf{X}}_N =\{ \underline{\bm{\tau}} \in \underline{\textbf{X}}: 
\textbf{n} \cdot \underline{\bm{ \tau} } = \bm{0}\quad \text{on} \quad \Gamma_N\},
\end{equation}
and define 
\begin{equation}
\bar{L}^2(\Omega) =
\begin{cases}
L^2(\Omega) &\quad \text{if} \quad\Gamma_N\neq \emptyset\\
\{ \bm{\varphi} \in L^2(\Omega) :
\int_{\Omega} \bm{\varphi}d\textbf{x} = 
0\} &\quad \text{if} \quad\Gamma_N= \emptyset\\
\end{cases}
\end{equation}
A least-squares formulation of Eq.~\eqref{eq:elasticity_with_vorticity} was considered in~\cite{Cai}.
By this approach, a three-field formulation seeks a minimizer of the functional
\begin{equation}
\mathcal{G}(\underline{\bm{\tau}},\textbf{v},\bm{\varphi};\textbf{f})=
||\mathcal{A}\underline{\bm{\tau}}-\nab{\bf{v}}
+
(-1)^d\underline{\bm{\chi}} \bm{\varphi} ||_0^2
+
||\textbf{div}\underline{\bm{\tau}}+\textbf{f}||_0^2
+
||\text{as}\underline{\bm{\tau}}||_0^2
\end{equation}
for all $(\underline{\bm{\tau}}, \bf{v},\bm{\varphi})$ $\in \underline{ \textbf{X} }_N \times H^1_{0,D}(\Omega)^d \times \bar{L}^2(\Omega)$.

%%%%%%%%%%%%%%%%%%%%%%%%%%%%%%%%%%%%%%%%%%%%%%%%%%%%%%%%%%%%%%%%%%%%%%%%%%%%%%%%
\section{Three-field eigenvalue formulation} 
\label{sec:three_field_eigenvalue_formulation}
%%%%%%%%%%%%%%%%%%%%%%%%%%%%%%%%%%%%%%%%%%%%%%%%%%%%%%%%%%%%%%%%%%%%%%%%%%%%%%%%
In this section we describe the three field formulation and its corresponding eigenvalue problem. The minimization of the functional $\mathcal{G}(\underline{\bm{\tau}},\bf{v},\bm{\varphi};\textbf{f})$ gives rise to the variational formulation: find $(\underline{\bm{\sigma}}, \textbf{u}, \bm{\psi}) \in \underline{ \textbf{X} }_N \times H^1_{0,D} (\Omega)^d \times\bar{L}^2(\Omega)$ such that
\begin{equation}\label{eq:variationa_form_source}
\begin{cases}
(\mathcal{A}\underline{\bm{\sigma}},\mathcal{A}\underline{\bm{\tau}})
+
(\textbf{div}\underline{\bm{\sigma}},\textbf{div}\underline{\bm{\tau}})
+
(\text{as}\underline{\bm{\sigma}},\text{as}\underline{\bm{\tau}})\\

\qquad -(\nab\textbf{u}, \mathcal{A}\underline{\bm{\tau}})
+
(-1)^d(\underline{\bm{\chi}}\bm{\psi}, \mathcal{A}\underline{\bm{\tau}})
=
-(\textbf{f},\textbf{div}\underline{\bm{\tau}})
&\quad \forall \underline{\bm{\tau}} \in \underline{\textbf{X}}_N \\
-
(\mathcal{A}\underline{\bm{\sigma}}, \nab \textbf{v})
+
(\nab\textbf{u},  \nab\textbf{v})
-
(-1)^d(\underline{\bm{\chi}}\bm{\psi}, \nab \textbf{v})
=
\bf{0}
&\quad \forall\textbf{v} \in H^{1}_{0,D}(\Omega)^d \\
(-1)^d(\mathcal{A}\underline{\bm{\sigma}},\underline{\bm{\chi}}\bm{\varphi})
-
(-1)^d(\nab\textbf{u}, \underline{\bm{\chi}}\bm{\varphi})
+
(\underline{\bm{\chi}}\bm{\psi},\underline{\bm{\chi}}\bm{\varphi})
=
\bf{0}
&\quad \forall \bm{\varphi} \in \bar{L}^2(\Omega) \\
\end{cases}
\end{equation}
Eigenvalue problems based on the finite element least-squares formulations was first studied for the Laplacian in~\cite{FB_eig} and then investigated for the linear elasticity problem in~\cite{FB_LE}.
We consider the spectrum of the solution operator associated with our formulation.
This means that we replace the source term $\bf{f}$ with $\gamma \textbf{u}$, where $\gamma$ is the eigenvalue.
The problem in this case reads: find $(\gamma,\textbf{u}) \in \mathbb{C}\times H_{0,D}^1(\Omega)^d$ such that $\textbf{u}\neq \textbf0$, and for some underline ${\bm{\sigma}}\in \underline{\textbf{X}}_N$ and $\bm{\psi}\in\bar{L}^2$ we have
\begin{equation}\label{eq:variationa_form_eigs}
\begin{cases}
({\mathcal{A}}\underline{\bm{\sigma}},\mathcal{A}\underline{\bm{\tau}})
+
(\textbf{div}\underline{\bm{\sigma}},\textbf{div}\underline{\bm{\tau}})
+
(\text{as}\underline{\bm{\sigma}},\text{as}\underline{\bm{\tau}})\\
\qquad-(\nab\textbf{u}, \mathcal{A}\underline{\bm{\tau}})
+
(-1)^d(\underline{\bm{\chi}}\bm{\psi}, \mathcal{A}\underline{\bm{\tau}})
=
-\gamma(\textbf{u},\bf{div}\underline{\bm{\tau}})
&\quad \forall \underline{\bm{\tau}} \in \underline{\textbf{X}}_N \\
-
(\mathcal{A}\underline{\bm{\sigma}}, \nab \textbf{v})
+
(\nab\textbf{u},  \nab\textbf{v})
-
(-1)^d(\underline{\bm{\chi}}\bm{\psi}, \nab \bf{v})
=
\bf{0}
&\quad \forall\textbf{v} \in H^{1}_{0,D}(\Omega)^d\\
(-1)^d(\mathcal{A}\underline{\bm{\sigma}},\underline{\bm{\chi}}\bm{\varphi})
-
(-1)^d(\nab\textbf{u}, \underline{\bm{\chi}}\bm{\varphi})
+
(\underline{\bm{\chi}}\bm{\psi},\underline{\bm{\chi}}\bm{\varphi})
=
\textbf{0}
&\quad \forall \bm{\varphi} \in \bar{L}^2(\Omega) \\
\end{cases}
\end{equation}
The above eigenvalue problem can be viewed in terms of operator in the following manner
\begin{equation}
\begin{cases}
A:(\mathbf{\mathcal{A}}\underline{\bm{\sigma}},\mathbf{\mathcal{A}}\underline{\bm{\tau}})
+
(\bf{div}\underline{\bm{\sigma}},\bf{div}\underline{\bm{\tau}})
+
(\text{as}\underline{\bm{\sigma}},\text{as}\underline{\bm{\tau}})\\
B:-(\mathbf{\mathcal{A}}\underline{\bm{\sigma}},
\nab\bf{v})\\
C:(-1)^d(\mathbf{\mathcal{A}}\underline{\bm{\sigma}},
\underline{\bm{\chi}}\bm{\varphi})\\
D:(\nab\bf{u},  \nab\bf{v}) \\
E:-(-1)^{d}(\nab\bf{u}, \underline{\bm{\chi}}\bm{\varphi})\\
F:(\underline{\bm{\chi}}\bm{\psi},\underline{\bm{\chi}}\bm{\varphi})\\
G:-(\bf{u},\bf{div}\underline{\bm{\tau}})
\end{cases}
\end{equation}
Thus, formulation \eqref{eq:variationa_form_eigs} involves a 3-by-3 block operators as follow
\begin{equation}\label{eq:matrix_structure}
\begin{pmatrix} 
A & B^t & C^t\\ 
B & D & E^t \\
C & E & F  
\end{pmatrix}
\begin{pmatrix} 
x \\ 
y \\
z
\end{pmatrix}
= \gamma
\begin{pmatrix} 
0 & G & 0\\
0 & 0 & 0\\
0 & 0 & 0
\end{pmatrix}
\begin{pmatrix} 
x \\ 
y \\
z
\end{pmatrix}.
\end{equation}
The system of equations \eqref{eq:matrix_structure} has clearly a symmetric left hand side.
On the contrary, the right hand side is not symmetric and singular.
Thus, the discrete spectrum is expected to contain complex and infinite eigenvalues. 
%%%%%%%%%%%%%%%%%%%%%%%%%%%%%%%%%%%%%%%%%%%%%%%%%%%%%%%%%%%%%%%%%%%%%%%%%%%%%%%%
\section{Finite element approximation}
\label{sec:FiniteElementApproximation}
%%%%%%%%%%%%%%%%%%%%%%%%%%%%%%%%%%%%%%%%%%%%%%%%%%%%%%%%%%%%%%%%%%%%%%%%%%%%%%%%
Let $\Omega$ be a polygonal domain, $h$ be the mesh-size and $\mathcal{T}_h=\{K\}$ be finite element partition of domain with elements being triangular.
The finite elements spaces used were proposed in \cite{Cai}. 
These finite element spaces are Raviart–Thomas space of degree $k$  ($RT_k$) for the stress, 
standard (conforming) continuous piecewise polynomials of degree $k + 1$ for the displacement denoted by $CG_{k+1}$, and finally discontinuous piecewise polynomials of degree $k$ for the rotation, denoted by $DG_k$.
These spaces are
$$\Sigma_h =\bigg \{ \underline{\bm{\tau}} \in \Sigma_k: \underline{\bm{\tau}}|_K \in {RT}_k(K)^d, \quad \forall K \in \mathcal{T}_h \bigg\} \subset \underline{\bm{X}}_N$$
$$U_h=\bigg \{ \textbf{v} \in \mathcal{C}^0(\Omega)^d : \textbf{v}|_K \in P_{(k+1)}(K)^d
,\quad \forall K \in \mathcal{T}_h, \textbf{v}=\bm{0} \text{ on } \Gamma_D \bigg\} 
\subset H_{D,0}^{1}(\Omega)^d$$
$$ \Phi_h = \bigg\{ \bm{\varphi} \in \bar{L}^2(\Omega):\bm{\varphi}|_K \in P_k(K), \quad\forall  K \in \mathcal{T}_h, \int_\Omega \bm{\varphi} dx = 0
\text{ if } \Gamma_N = \emptyset \bigg\} \subset L^2_D(\Omega)$$
Given those finite dimensional subspaces, the Galerkin approximation of Eq.~\eqref{eq:variationa_form_eigs} is then: find $(\gamma_h, \textbf{u}_h) \in\mathbb{C} \times U_h$ $\textbf{u}_h\ne0$ such that for some $\underline{\bm{\sigma}}_h\in \Sigma_h$ and some $\bm{\psi}_h \in \Phi_h$ we have
\begin{equation}\label{eq:discrete_variationa_form_eigs}
\begin{cases}
(\mathcal{A}\underline{\bm{\sigma}}_h,\mathcal{A}\underline{\bm{\tau}})
+
(\textbf{div}\underline{\bm{\sigma}}_h,\textbf{div}\underline{\bm{\tau}})
+
(\text{as}\underline{\bm{\sigma}}_h,\text{as}\underline{\bm{\tau}})\\

\qquad -(\nab\textbf{u}_h, \mathcal{A}\underline{\bm{\tau}})
+
(-1)^d(\underline{\bm{\chi}}\bm{\psi}_h, \mathbf{\mathcal{A}}\underline{\bm{\tau}})
=
-\gamma_h(\textbf{u}_h,\textbf{div}\underline{\bm{\tau}})
&\quad \forall \underline{\bm{\tau}} \in \Sigma_h \\
-
(\mathcal{A}\underline{\bm{\sigma}}_h, \nab \textbf{v})
+
(\nab\textbf{u}_h,  \nab\textbf{v})
-
(-1)^d(\underline{\bm{\chi}}\bm{\psi}_h, \nab \bf{v})
=
\bf{0}
&\quad \forall\textbf{v} \in U_h \\
(-1)^d(\mathcal{A}\underline{\bm{\sigma}}_h,\underline{\bm{\chi}}\bm{\varphi})
-
(-1)^d(\nab\textbf{u}_h, \underline{\bm{\chi}}\bm{\varphi})
+
(\underline{\bm{\chi}}\bm{\psi}_h,\underline{\bm{\chi}}\bm{\varphi})
=
0
&\quad \forall \bm{\varphi} \in \Phi_h \\
\end{cases}
\end{equation}
The structural formulation of the above algebraic system is exactly the one that corresponds to 
\eqref{eq:matrix_structure}. 
This has the form
\begin{equation}
M v =\lambda N v,
\end{equation}
which is a generalized eigenvalue problem.
The vector $(x,y,z)^t$ in Eq.~\eqref{eq:matrix_structure} corresponds to $(\bm{\hat{\sigma}}_h, \hat{\textbf{u}}_h,\hat{\bm{\psi}}_h)^t$.
In~\cite{FB_eig} and~\cite{FB_LE} the interested reader can have some information on how to solve those problems, even if our main interest is not on the solution of the algebraic system.
The characterization of the computed eigenvalues and the strategy to solve such a system (with a singular right hand side matrix) are also discussed.

In the case of the three field formulation, by looking at the matrix structure~\ref{eq:matrix_structure}, it is possible to rewrite the system in the following manner:
\begin{equation}\label{sys:3_variables_Algebric_linear_system}
\begin{cases}
A\bm{\hat{\sigma}}_h + B^t \hat{\textbf{u}}_h + C^t \hat{\bm{\psi}}_h= \gamma_h G \hat{\bf{u}}_h\\
B\hat{\bm{\sigma}}_h + D \hat{\textbf{u}}_h + E^t \hat{\bm{\psi}}_h = 0\\
C\hat{\bm{\sigma}}_h + E \hat{\textbf{u}}_h + F \hat{\bm{\psi}}_h = 0
\end{cases} 
\end{equation}
Using the last equation in System~\eqref{sys:3_variables_Algebric_linear_system}, with $F$ being invertible, gives
\begin{equation}
\hat{\bm{\psi}}_h = - F^{-1} ( C \hat{\bm{\sigma}}_h + E \hat{\bf{u}}_h).
\end{equation}

Substituting the value of $\hat{\bm{\psi}}_h$ in \eqref{sys:3_variables_Algebric_linear_system} we deduce the following algebraic linear system
\begin{equation}\label{sys:2_variables_Algebric_linear_system}
\begin{cases}
\gA\bm{\hat{\sigma}}_h + \gB \hat{\bf{u}}_h = \gamma_h G \hat{\bf{u}}_h\\
\gC\hat{\bm{\sigma}}_h + \gD \hat{\bf{u}}_h = 0\\
\end{cases} 
\end{equation}
where the expression of the involved matrices is given by
\begin{equation}
\begin{cases}
\gA = A - C^t F^{-1}C \\
\gB = B^t-C^t F^{-1}E \\
\gC = B - E^t F^{-1}C \\
\gD = D - E^t F^{-1}E
\end{cases}
\end{equation}
By looking into the Schur complement for the non-vanishing displacement, and since the matrix $A$ is non-singular, this give rise to the subsequent expression
\begin{equation}\label{eq:schur_threefield_displacment}
(\underbrace{\gC \text{ } \gA^{-1}\gB - \gD}_{M}) \hat{\textbf{u}}_h = 
\gamma_h \underbrace{\gC \text{ }\gA^{-1}G}_N \hat{\bf{u}}_h,
\end{equation}
which has the structure of a generalized eigenvalue problem.
%%%%%%%%%%%%%%%%%%%%%%%%%%%%%%%%%%%%%%%%%%%%%%%%%%%%%%%%%%%%%%%%%%%%%%%%%%%%%%%%
\section{Numerical Results}
\label{se:num}
%%%%%%%%%%%%%%%%%%%%%%%%%%%%%%%%%%%%%%%%%%%%%%%%%%%%%%%%%%%%%%%%%%%%%%%%%%%%%%%%
The aim of the paper is to present some numerical results that confirm the theoretical analysis on the convergence of the spectrum discussed in~\cite{FB_LE}.
We continue the investigation on the approximated eigenvalues for the three-field formulation when specific meshes are selected.
A similar investigation was conducted in~\cite{LZ_DB} concerning the two-field formulation proposed in~\cite{FB_LE}.
The general theory about the spectrum of operators derived from least-squares was originally studied in~\cite{FB_eig}.

The linear elasticity problem investigated in this article is associated with the compliance tensor defined in~\eqref{eq:compliance_tensor}.
As explained previously, we are concerned with the situation where the material tends to approach the incompressible limit.
We study the problem in $\mathbb{R}^2$ with the Lam\'e constant $\mu$ being $1$.
Hence, the material being compressible/incompressible is associated with the value of $\lambda$; we start by taking $\lambda=1$. We then gradually move to the incompressible limit by choosing the values $\lambda=10^r$ with $r=2,4,8$.

Two domains are examined, a square and an L-shaped domain.
Each domain is associated with three mesh sequences.
The studied domains are subdivided into a finite number of elements on each side denoted by $N$. The larger the value of $N$ the finer the mesh.

The chosen spaces are the finite elements defined in Section \ref{sec:FiniteElementApproximation} with $k=0$.
These spaces are $RT_0$ for the stress, $CG_1$ for the displacement and $DG_0$ for the rotation. In order to construct the spaces and build the matrices defined in~\eqref{eq:matrix_structure}, the FEniCS project~\cite{Fenics} was utilized.
The Schur complement of the displacement resulted in~\ref{eq:schur_threefield_displacment} is then extracted in order to solve the eigenvalue problem.

In what follows we first present the rate of convergence for the first eigenvalue.
The reference solutions used in estimating the rate for both domains were provided in~\cite{LZ_DB} (see Table 1,2).
The behavior of the spectrum in the complex plane is then studied with different scales of the
axes.
This is done to better describe the spread of eigenvalues for different values of $\lambda$.
Since we didn't observe significant differences between the cases $\lambda=10^4 $ and $\lambda=10^8$, we only illustrate the results for $\lambda=10^8$.
Each domain is reported separately and we finally conclude with a comparison between all cases.
we note that with abuse of notation, we call the modules of the complex part being close/far from the real axes in short by small/large imaginary parts.
%%%%%%%%%%%%%%%%%%%%%%%%%%%%%%%%%%%%%%%%%%%%%%%%%%%%%%%%%%%%%%%%%%%%%%%%%%%%%%%%
\subsection{Numerical results on the square}
\label{sub:numerical_results_on_the_square}
\begin{figure}[ht!]
\centering
\begin{subfigure}[b]{0.32\textwidth}
\centering
\includegraphics[width=\textwidth]{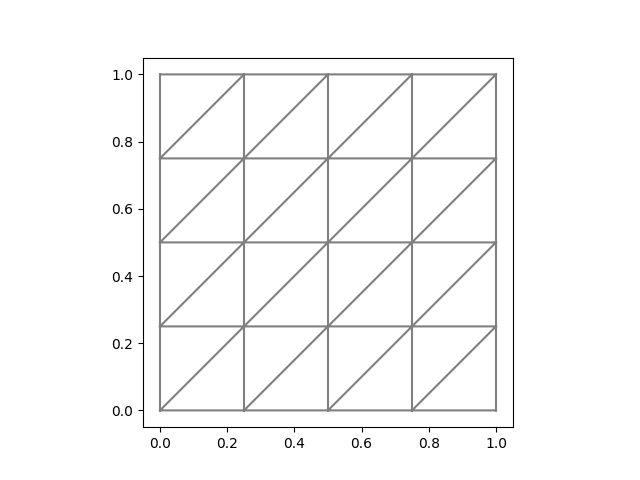}
\caption{Right}
% \label{fig:right_square}
\end{subfigure}
\hfill
\begin{subfigure}[b]{0.32\textwidth}
\centering
\includegraphics[width=\textwidth]{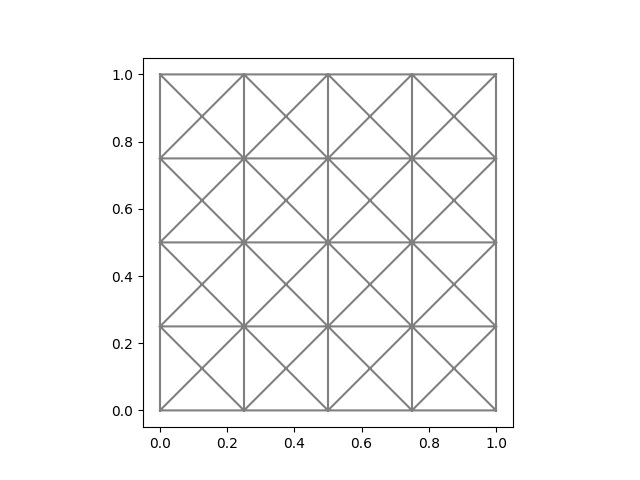}
\caption{Crossed}
% \label{fig:crossed_square}
\end{subfigure}
\hfill
\begin{subfigure}[b]{0.32\textwidth}
\centering
\includegraphics[width=\textwidth]{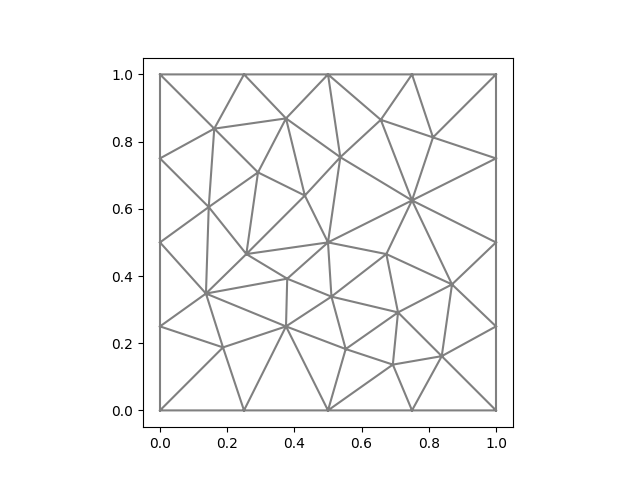}
\caption{Nonuniform}
% \label{fig:nonuniform_square}
\end{subfigure}
\caption{Meshes for the unit square domain with $N=4$}
\label{fig:square_meshes}
\end{figure}

%%%%%%%%%%%%%%%%%%%%%%%%%%%%%%%%%%%%%%%%%%%%%%%%%%%%%%%%%%%%%%%%%%%%%%%%%%%%%%%%
We start our investigation with the square domain $\Omega=]0,1[^2$.
Three kinds of meshes are examined: Right, Crossed and Nonuniform with their structure being non-symmetric uniform, symmetric uniform and non-symmetric non-uniform, respectively.
An example of such meshes with $N=4$ are plotted in  Fig.~\ref{fig:square_meshes}. 
%%%
\begin{figure}[ht!]
\centering
\begin{subfigure}[b]{0.32\textwidth}
\centering
\includegraphics[width=\textwidth]{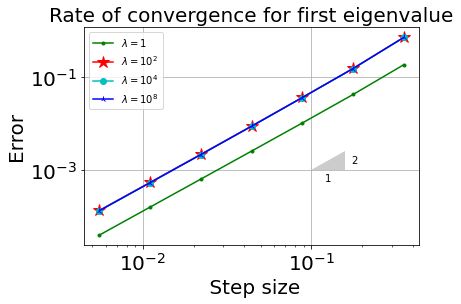}
\caption{Right}
% \label{fig:right_square}
\end{subfigure}
\hfill
\begin{subfigure}[b]{0.32\textwidth}
\centering
\includegraphics[width=\textwidth]{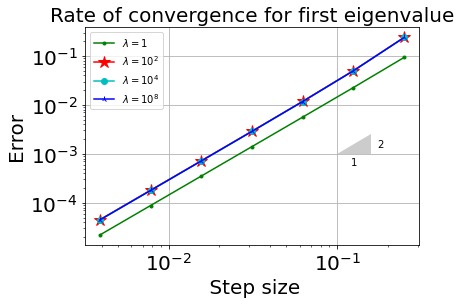}
\caption{Crossed}
% \label{fig:crossed_square}
\end{subfigure}
\hfill
\begin{subfigure}[b]{0.32\textwidth}
\centering
\includegraphics[width=\textwidth]{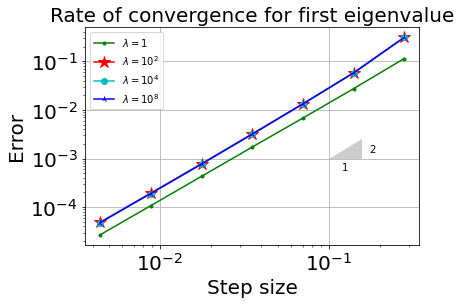}
\caption{Nonuniform}
% \label{fig:nonuniform_square}
\end{subfigure}
\caption{Rate of convergence for the first eigenvalue on a square}
\label{fig:Rate_of_convergence_for_the_first_eigenvalue_on_a_square}
\end{figure}
When using the spaces specified previously, authors in~\cite{FB_LE} proved that the rate of convergence is of second order.
Figure~\ref{fig:Rate_of_convergence_for_the_first_eigenvalue_on_a_square} reports the rate of the first eigenvalue approximated which is aligned with the proven theory.
\begin{figure}[ht!]
\centering
\begin{subfigure}[b]{0.3\textwidth}
\centering
\includegraphics[width=\textwidth]{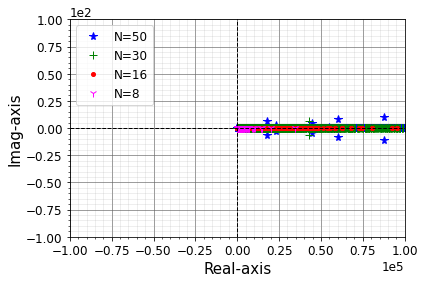}
\caption{$\lambda = 1$}
\label{fig:lam1_right_refine}
\end{subfigure}
\hfill
\begin{subfigure}[b]{0.3\textwidth}
\centering
\includegraphics[width=\textwidth]{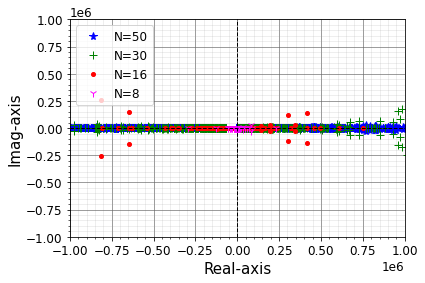}
\caption{$\lambda = 100$}
\label{fig:lam100_right_refine}
\end{subfigure}
\hfill
\begin{subfigure}[b]{0.3\textwidth}
\centering
\includegraphics[width=\textwidth]{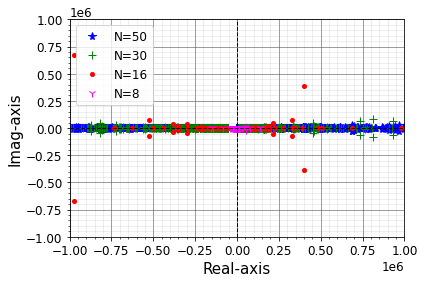}
\caption{$\lambda = 10^8$}
\label{fig:lam10pow8_right_refine}
\end{subfigure}
\caption{Spread of Eigenvalues for Right mesh on unit square} 
\label{fig:right_square_refine}
\end{figure}
We now discuss the distribution of the spectrum in the complex plane for the square, emphasizing the fact that the exact eigenvalues are real and positive.
Starting with the Right mesh and the material being solid elastic, Fig.~\ref{fig:lam1_right_refine} shows that the eigenvalues are concentrated to the right half of the plane, having positive real parts and small complex values appear as we refine.
On the contrary, for larger values of $\lambda$, as Figs.~\ref{fig:lam100_right_refine} and~\ref{fig:lam10pow8_right_refine} show, negative real eigenvalues appear.
As the mesh is refined, the eigenvalues either converge to positive real numbers or diverge with growing modulus.
\begin{figure}[ht!]
\centering
\begin{subfigure}[b]{0.3\textwidth}
\centering
\includegraphics[width=\textwidth]{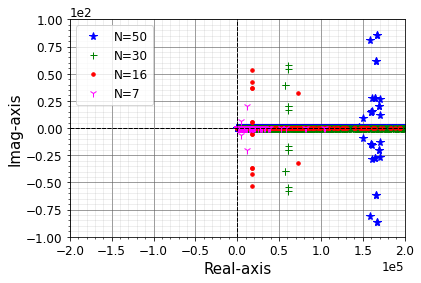}
\caption{$\lambda = 1$}
\label{fig:lam1_crossed_refine}
\end{subfigure}
\hfill
\begin{subfigure}[b]{0.3\textwidth}
\centering
\includegraphics[width=\textwidth]{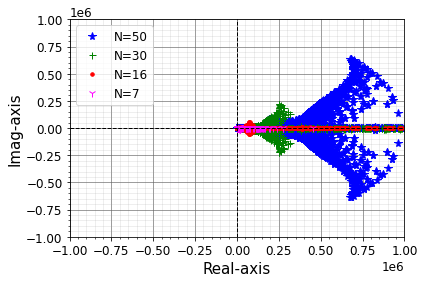}
\caption{$\lambda = 100$}
\label{fig:lam100_crossed_refine}
\end{subfigure}
\hfill
\begin{subfigure}[b]{0.3\textwidth}
\centering
\includegraphics[width=\textwidth]{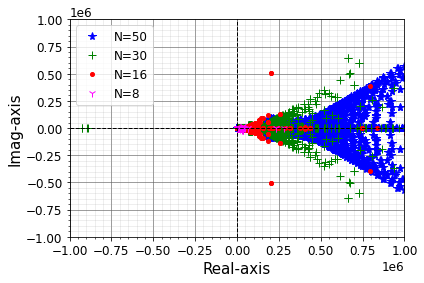}
\caption{$\lambda = 10^8$}
\label{fig:lam10pow8_crossed_refine}
\end{subfigure}
\caption{Spread of Eigenvalues for Crossed mesh on unit square} 
\label{fig:crossed_square_refine}
\end{figure}

Figure~\ref{fig:crossed_square_refine} illustrates the distribution of eigenvalues for the Crossed mesh.
In this case, the symmetry of the mesh plays a crucial role in the layout of eigenvalues.
All eigenvalues are concentrated to the right half of the complex plane.
We start with Fig.~\ref{fig:lam1_crossed_refine} which presents the distribution of the eigenvalues for $\lambda=1$.
It is evident that with a crossed mesh more eigenvalues with large complex part appear than what happens with the Right meshes.
As the Lam\'e constant $\lambda$ becomes larger (cf. Figs.~\ref{fig:lam100_crossed_refine} and~\ref{fig:lam10pow8_crossed_refine}) the eigenvalues spread more and move far away from the origin.
Also in this case the result confirm the convergence as $h$ is refined.

\begin{figure}[ht!]
\centering
\begin{subfigure}[b]{0.3\textwidth}
\centering
\includegraphics[width=\textwidth]{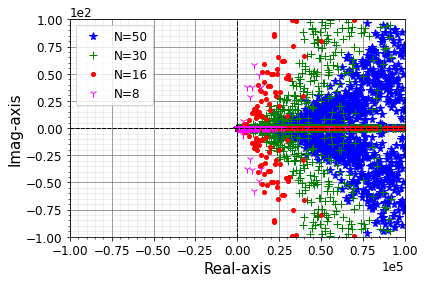}
\caption{$\lambda = 1$}
\label{fig:lam1_polygon_refine}
\end{subfigure}
\hfill
\begin{subfigure}[b]{0.3\textwidth}
\centering
\includegraphics[width=\textwidth]{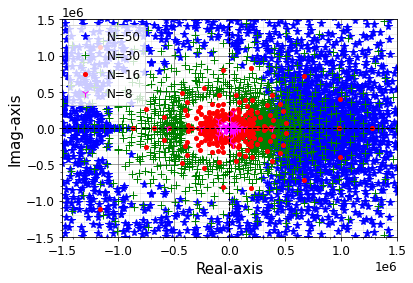}
\caption{$\lambda = 100$}
\label{fig:lam100_polygon_refine}
\end{subfigure}
\hfill
\begin{subfigure}[b]{0.3\textwidth}
\centering
\includegraphics[width=\textwidth]{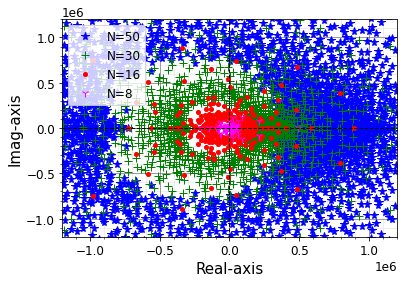}
\caption{$\lambda = 10^8$}
\label{fig:lam10pow8_polygon_refine}
\end{subfigure}
\caption{Spread of Eigenvalues for Nonuniform mesh on unit square} 
\label{fig:polygon_square_refine}
\end{figure}

In the case of Nonuniform meshes, different behaviors are observed.
This is shown in Fig~\ref{fig:polygon_square_refine}.
The spectrum in this case is more disturbed and scattered allover the complex plane.
It is only the case when a small value of $\lambda$ is considered where the material is compressible, the eigenvalues appear to the right half of the complex plane with positive real values.
For larger values of $\lambda$ that is $100$ or $10^8$, the spectrum gets disperse with the eigenvalues being more concentrated to the right half of the complex plane.

%%%%%%%%%%%%%%%%%%%%%%%%%%%%%%%%%%%%%%%%%%%%%%%%%%%%%%%%%%%%%%%%%%%%%%%%%%%%%%%%
\subsection{Numerical results on the L-shaped domain} 
\label{sub:Numerical_results_on_the_L-shaped_domain}
%%%%%%%%%%%%%%%%%%%%%%%%%%%%%%%%%%%%%%%%%%%%%%%%%%%%%%%%%%%%%%%%%%%%%%%%%%%%%%%%

\begin{figure}[ht!]
\centering
\begin{subfigure}[b]{0.3\textwidth}
\centering
\includegraphics[width=\textwidth]{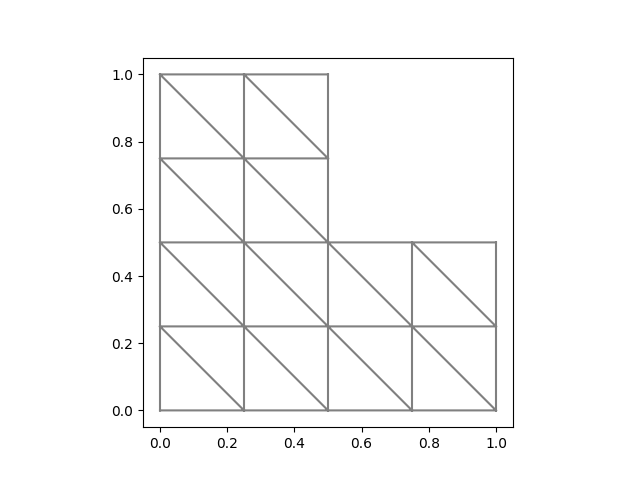}
\caption{Left}
\label{fig:lshape_left}
\end{subfigure}
\hfill
\centering
\begin{subfigure}[b]{0.3\textwidth}
\centering
\includegraphics[width=\textwidth]{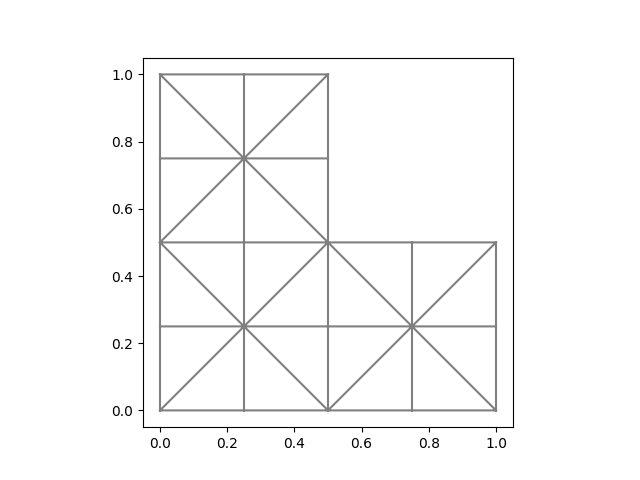}
\caption{Uniform}
\label{fig:lshaped_uniform}
\end{subfigure}
\hfill
\begin{subfigure}[b]{0.3\textwidth}
\centering
\includegraphics[width=\textwidth]{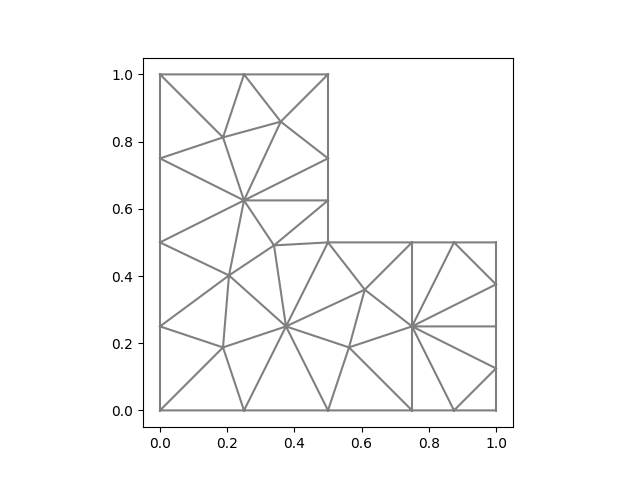}
\caption{Nonuniform}
\label{fig:lshaped_Non_uniform}
\end{subfigure}
\caption{Meshes for the L-shaped domain with $N=4$} 
\label{fig:lshape_mesh}
\end{figure}
We also inspect and study the L-shaped domain where $\Omega = ]0, 1[^2  \setminus ([0.5, 1[ \times [0.5, 1[)$.
Also in this case we consider three different mesh sequences which are Left, Uniform and Nonuniform meshes.
These are illustrated in  Fig.~\ref{fig:lshape_mesh} with $N=4$.
It is well known that the eigenfunctions computed on the L-shaped domain may have singularities due to the re-entered corner.
We are considering the first eigenmode which is known to correspond to a singular solution.
\begin{figure}[ht!]
\centering
\begin{subfigure}[b]{0.32\textwidth}
\centering
\includegraphics[width=\textwidth]{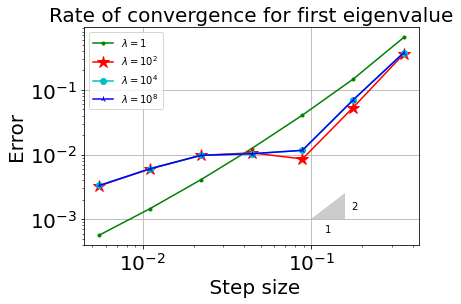}
\caption{Left}
\label{fig:lshape_left_rate}
\end{subfigure}
\hfill
\begin{subfigure}[b]{0.32\textwidth}
\centering
\includegraphics[width=\textwidth]{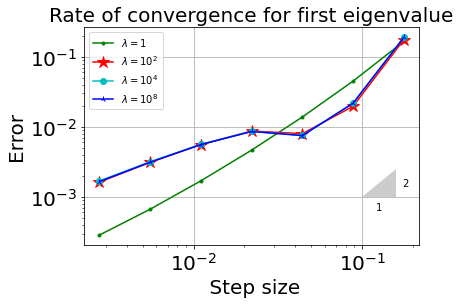}
\caption{Uniform}
\label{fig:lshape_uniform_rate}
\end{subfigure}
\hfill
\begin{subfigure}[b]{0.32\textwidth}
\centering
\includegraphics[width=\textwidth]{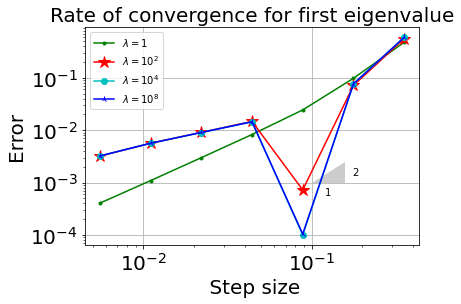}
\caption{Nonuniform}
\label{fig:lshape_nonuniform_rate}
\end{subfigure}
\caption{Rate of convergence for the first eigenvalue on the L-shaped domain}
\label{fig:Rate_of_convergence_for_the_first_eigenvalue_on_lshape}
\end{figure}
%%%
\begin{figure}
\centering
\begin{subfigure}[b]{0.45\textwidth}
\centering
\includegraphics[width=\textwidth]{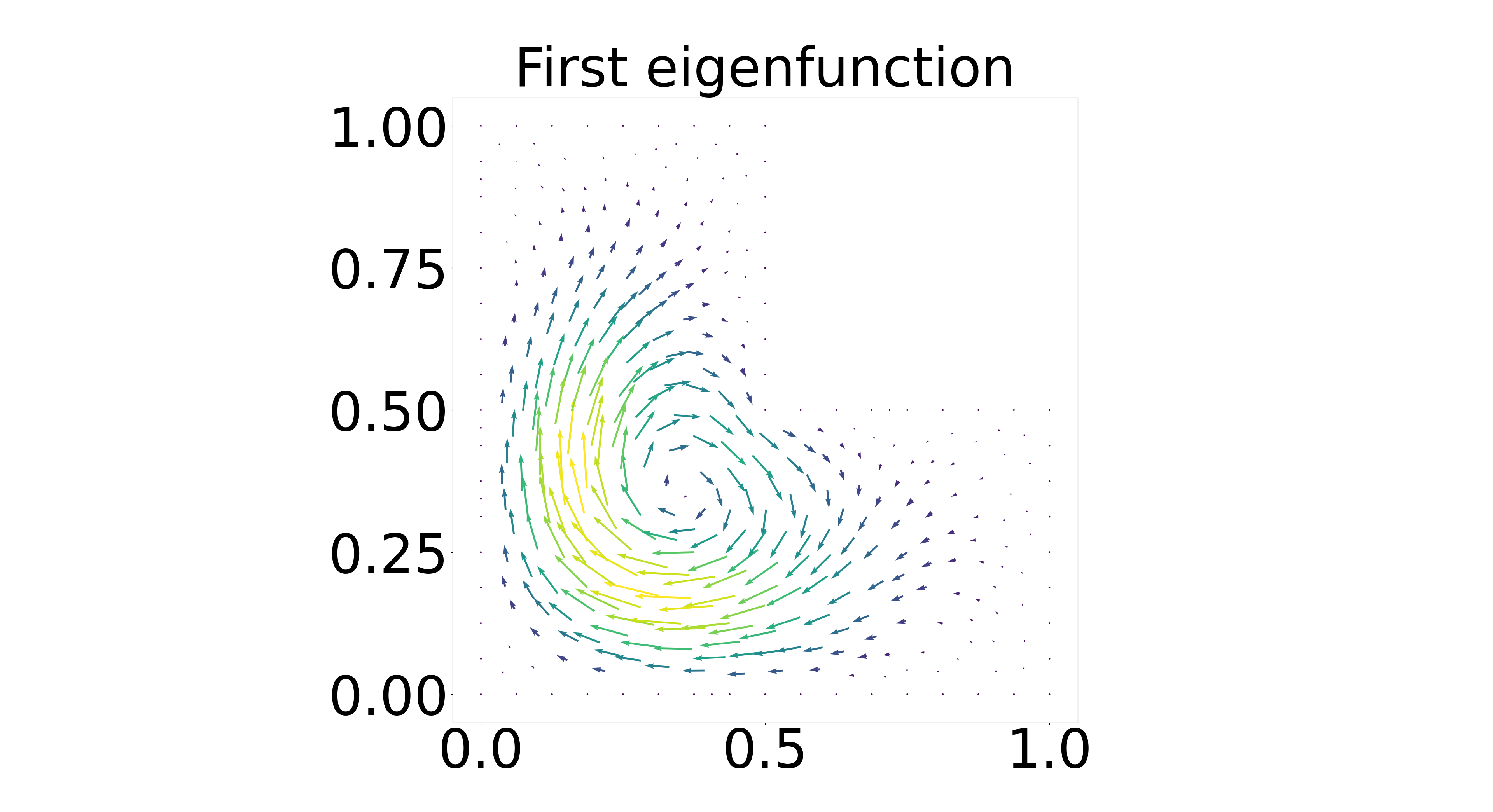}
\caption{$N=16$}
\label{fig:first_eigfun_nonuiform_lshape_N16_lam100}
\end{subfigure}
\begin{subfigure}[b]{0.45\textwidth}
\centering
\includegraphics[width=\textwidth]{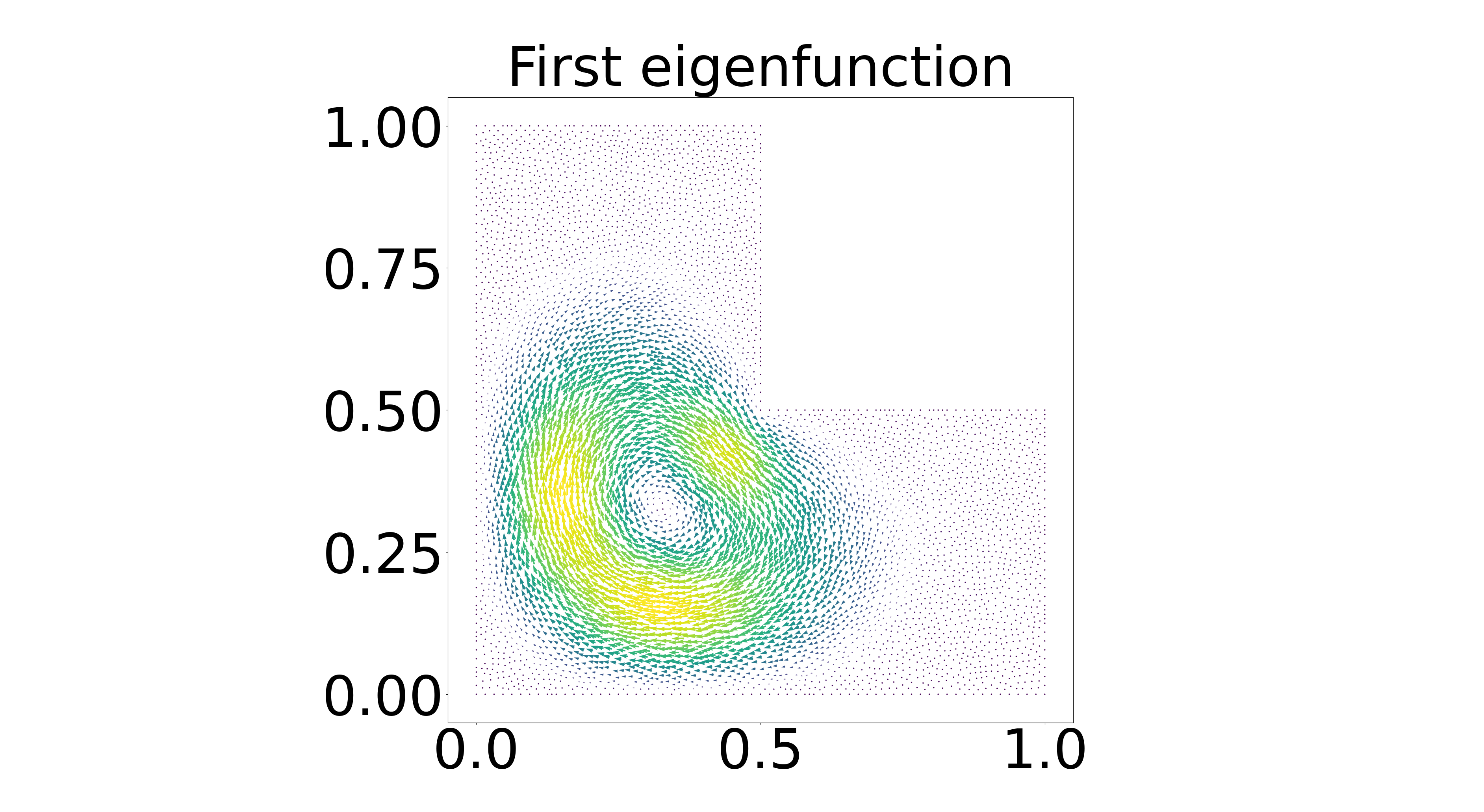}
\caption{$N=64$}
\label{fig:first_eigfun_nonuiform_lshape_N64_lam100}
\end{subfigure}
\caption{First eigenfunction for Nonuniform mesh with $\lambda=100$}
\label{fig:first_eigfun_nonuiform_lshape_lam100}
\end{figure}

Figure~\ref{fig:Rate_of_convergence_for_the_first_eigenvalue_on_lshape} presents the rate for three mesh structures defined on the L-shaped domain with different $\lambda$ values.
When computing the rate of convergence for the first eigenvalue of our problem, we observe the same pre-asymptotic phenomena (the convergence curve is not straight) as in~\cite{LZ_DB}.
This appears when $\lambda$ grows moving to the incompressible limit.
Thus, we choose to compute the first eigenfunction (Fig.~\ref{fig:first_eigfun_nonuiform_lshape_lam100}) on a Nonuniform mesh with $N=16$ and $64$ for $\lambda=100$ to observer what does the eigenfunction looks like.
Taking the eigenfunction for $N=16$ as Fig.~\ref{fig:first_eigfun_nonuiform_lshape_N16_lam100} shows, we clearly see that a coarse mesh does not represent the correct re-circulation of the vortex expected when moving to the incompressible case.
However, taking a finer mesh with $N=64$ in Fig.~\ref{fig:first_eigfun_nonuiform_lshape_N64_lam100} gives nearly the expected vorticity which in return gives the expected rate, being less that $2$, as seen in Figs.~\ref{fig:lshape_left_rate} and~\ref{fig:lshape_uniform_rate}.
Moreover, in Fig.~\ref{fig:lshape_nonuniform_rate} we notice that the pre-asymptotic behavior acts as super convergence between the second and third iteration.
We attribute this behavior to the closeness of the exact and approximated solutions of our used “accurate” eigenvalue which was computed by the \texttt{\uppercase{Solve-Estimate-Mark-Refine}} strategy (cf.~\cite{LZ_DB} Sec 5). 
\begin{figure}[ht!]
\centering
\begin{subfigure}[b]{0.32\textwidth}
\centering
\includegraphics[width=\textwidth]{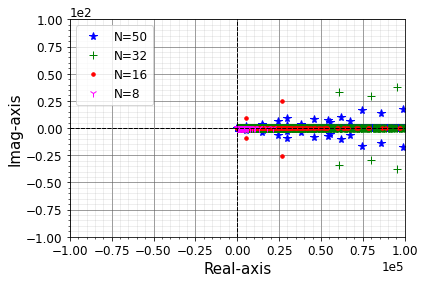}
\caption{$\lambda =1$}
\label{fig:lshape_left_lam1}
\end{subfigure}
\hfill
\begin{subfigure}[b]{0.32\textwidth}
\centering
\includegraphics[width=\textwidth]{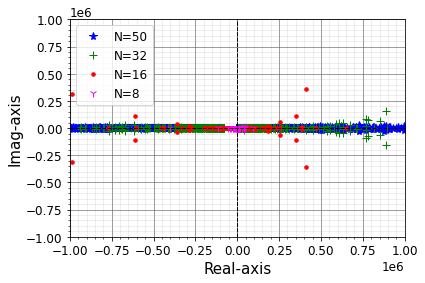}
\caption{$\lambda=100$}
\label{fig:lshape_left_lam100}
\end{subfigure}
\hfill
\begin{subfigure}[b]{0.32\textwidth}
\centering
\includegraphics[width=\textwidth]{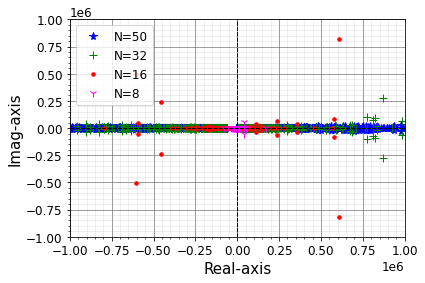}
\caption{$\lambda =10^8$}
\label{fig:lshape_left_lam10pow8}
\end{subfigure}
\caption{Spread of Eigenvalues for Left mesh on L-shaped domain}
\label{fig:Spread_of_Eigenvalues_for_Left_mesh_on_L-shape_domain}
\end{figure}

We now consider the distribution of eigenvalues for the L-shaped domain.
Figure~\ref{fig:Spread_of_Eigenvalues_for_Left_mesh_on_L-shape_domain} illustrates the spread when considering the Left mesh structure.
We see in such a case that the behavior of the eigenvalues is the same as in the Right mesh on the square. When $\lambda=1$ more complex eigenvalues appear when comparing it to the Right mesh on a square with its real part being positive as Fig.~\ref{fig:lshape_left_lam1} shows.
In general, the higher the $\lambda$ the more eigenvalues appear on both sides of the complex plane with small imaginary part.
\begin{figure}[ht!]
\centering
\begin{subfigure}[b]{0.32\textwidth}
\centering
\includegraphics[width=\textwidth]{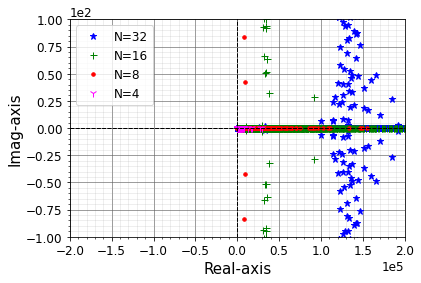}
\caption{$\lambda =1$}
\label{fig:lshape_uniform_lam1}
\end{subfigure}
\hfill
\begin{subfigure}[b]{0.32\textwidth}
\centering
\includegraphics[width=\textwidth]{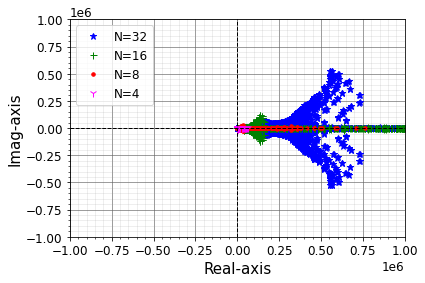}
\caption{$\lambda=100$}
\label{fig:lshape_uniform_lam100}
\end{subfigure}
\hfill
\begin{subfigure}[b]{0.32\textwidth}
\centering
\includegraphics[width=\textwidth]{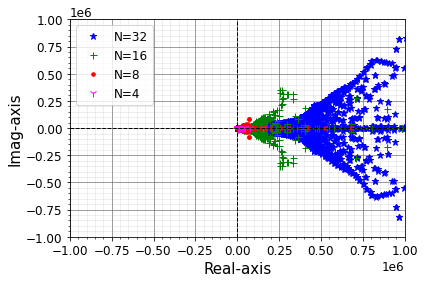}
\caption{$\lambda =10^8$}
\label{fig:lshape_uniform_lam10pow8}
\end{subfigure}
\caption{Spread of Eigenvalues for Uniform mesh on L-shaped domain}
\label{fig:Spread_of_Eigenvalues_for_Uniform_mesh_on_L-shape_domain}
\end{figure}

For the Uniform mesh case as Fig.~\ref{fig:Spread_of_Eigenvalues_for_Uniform_mesh_on_L-shape_domain}
shows, again the structure of the mesh plays a significant role in the distribution of the eigenvalues.
All eigenvalues have positive real parts. As $\lambda$ increases in this case, the spreading is more noticeable when comparing it to the crossed mesh.
Taking $N$ to be around $30$ as an example, we see that the eigenvalues in Fig.~\ref{fig:lshape_uniform_lam10pow8} in blue are more diffused when comparing it to the green eigenvalues in Fig.~\ref{fig:lam10pow8_crossed_refine} with the same axes for both figures.
Thus, the domain chosen and the mesh in hand plays a crucial role in the spread of eigenvalues with the method considered.
\begin{figure}[ht!]
\centering
\begin{subfigure}[b]{0.32\textwidth}
\centering
\includegraphics[width=\textwidth]{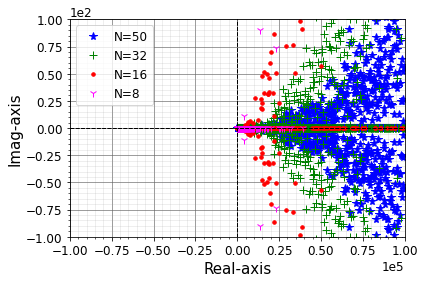}
\caption{$\lambda =1$}
\label{fig:lshape_Nonuniform_lam1}
\end{subfigure}
\hfill
\begin{subfigure}[b]{0.32\textwidth}
\centering
\includegraphics[width=\textwidth]{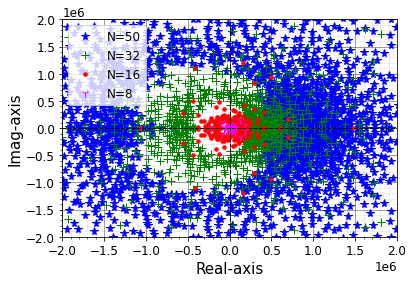}
\caption{$\lambda=100$}
\label{fig:lshape_Nonuniform_lam100}
\end{subfigure}
\hfill
\begin{subfigure}[b]{0.32\textwidth}
\centering
\includegraphics[width=\textwidth]{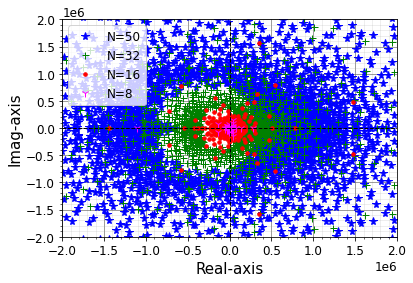}
\caption{$\lambda =10^8$}
\label{fig:lshape_Nonuniform_lam10pow8}
\end{subfigure}
\caption{Spread of Eigenvalues for Nonuniform mesh on L-shaped domain}
\label{fig:Spread_of_Eigenvalues_for_Nonuniform_mesh_on_L-shape_domain}
\end{figure}
\begin{figure}[ht!]
\centering
\begin{subfigure}[b]{0.45\textwidth}
\centering
\includegraphics[width=\textwidth]{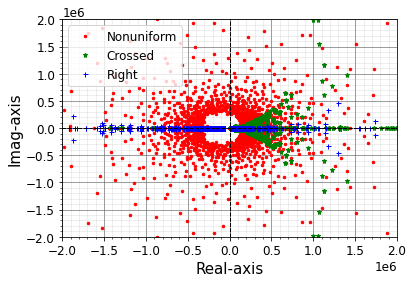}
\caption{Unit square with $N=30$}
\label{fig:allmeshes_unit_square_lam10pow8}
\end{subfigure}
\hfill
\begin{subfigure}[b]{0.45\textwidth}
\centering
\includegraphics[width=\textwidth]{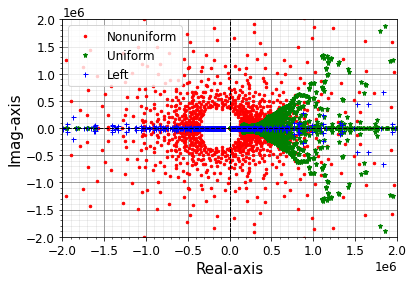}
\caption{L-shaped with $N=32$}
\label{fig:allmeshes_lshape_lam10pow8}
\end{subfigure}
\caption{Spread of Eigenvalues for all meshes}
\label{fig:Spread_of_Eigenvalues_for_all_meshes_lam10pow8}
\end{figure}

The Nonuniform mesh is the last case shown in Fig.~\ref{fig:Spread_of_Eigenvalues_for_Nonuniform_mesh_on_L-shape_domain}.
We observe the same layout of eigenvalues as in the case of Nonuniform mesh on the square.
All eigenvalues in this case are being spread except for $\lambda=1$.
In this case, the eigenvalues have positive real part and growing large imaginary as we refine.

For the sake of comparison between meshes and domains, we plot in Fig.~\ref{fig:Spread_of_Eigenvalues_for_all_meshes_lam10pow8} the distribution of the eigenvalues for the different meshes. We consider both domains when the material tends to the  incompressible case with $\lambda =10^8$.
Both figures show how the eigenvalues spread depending on the mesh chosen.
It is clear that in the case of the L-shaped domain (cf. Fig.~\ref{fig:allmeshes_lshape_lam10pow8}), the complex eigenvalues appear more scattered when comparing it to the unit square.

\bibliographystyle{plain}

\bibliography{qhe}

\end{document}